\documentclass[12pt]{amsart}
\usepackage{amsfonts,amsmath,amssymb,amsthm,comment,float,hyperref,mathrsfs,mathtools,multirow,times, color}

\newtheorem{theorem}{Theorem}[section]

\newtheorem{cor}[theorem]{Corollary}

\theoremstyle{definition}
\newtheorem{definition}[theorem]{Definition}

\theoremstyle{remark}
\numberwithin{equation}{section}

\allowdisplaybreaks

\newcommand\mylabel[1]{\label{#1}}
\newcommand\mybibitem[1]{\bibitem{#1}}
\newcommand\thm[1]{\ref{thm:#1}}

\newcommand\corol[1]{\ref{cor:#1}}

\newcommand\eqn[1]{(\ref{eq:#1})}
\newcommand\sect[1]{\ref{sec:#1}}

\newcommand\numberthis{\addtocounter{equation}{1}\tag{\theequation}}

\begin{document}

\title{Overpartitions with repeated smallest non-overlined part}

\author{Amita Malik}
\address{\newline \newline Department of Mathematics, The Pennsylvania State University, McAllister Building, University Park, PA 16802, USA \newline}
\email{malik@psu.edu}
\author{Rishabh Sarma}

\email{rishabh.sarma@psu.edu}

\maketitle

\begin{abstract}
Inspired by Andrews' and Bachraoui's work on partitions with repeated smallest part, we extend the concept to overpartitions. We study overpartitions with the restriction that the smallest non-overlined part appears exactly $k$ times and every overlined part is bigger than this part. We prove results expressing the generating functions of these overpartitions (and their subclass where no part has the same parity as the smallest part, among others) as linear combinations of the $q$-Pochhammer symbols with rational functions in $q$ as coefficients.
\end{abstract}

\section{Introduction}
\mylabel{section1}

In recent work \cite{An-Ba25}, Andrews and Bachraoui consider integer partitions whose smallest part is repeated exactly $k$-times and the remaining parts are distinct, denoted by $\mathrm{spt}k_d(n)$, and prove that their generating functions are linear combinations of the $q$-Pochhammer symbols with polynomials in $q$ as coefficients. 
\\\\
Let $N,n$ be non-negative integers. Throughout, we use the following standard notations.
\begin{align*}
    (a)_N = (a;q)_N &:= \prod_{k=0}^{N-1}(1-aq^k),\\
    (a)_{\infty} = (a;q)_{\infty} &:= \lim_{N\rightarrow \infty}(a)_N\,\,\text{where}\,\,\lvert q\rvert<1.
\end{align*}
\\
The following two results of Andrews and Bachraoui concern the class of partitions mentioned earlier.

\begin{theorem}
\mylabel{thm:An-Ba1}{\cite[Theorem 1]{An-Ba25}}
For any positive integer $k$, we have
\begin{align*}
\sum\limits_{n=1}^{\infty}\mathrm{spt}k_d(n)q^n=P_k(q)(-q;q)_{\infty}+(-1)^k(q;q)_{k-1},
\end{align*}
where 
\begin{align*}
P_k(q)=\begin{cases}
1 & \text{if } k=1, \\
(q^{k-1}-1) P_{k-1}(q)+q^{k-1}& \text{if } k>1.
\end{cases}
\end{align*}
\end{theorem}
\noindent
Moreover, a closed formula for $P_k(q)$ is given by the following theorem.
\begin{theorem}{\cite[Theorem 2]{An-Ba25}}
\mylabel{thm:An-Ba2}
For any positive integer $k$, we have
\begin{align*}
P_{2k-1}(q)=1+2\sum_{j=1}^{k-1}(-1)^jq^{j^2}(q^{2k-2j};q^2)_j,
\\
P_{2k}(q)=-1-2\sum_{j=1}^{k}(-1)^jq^{j^2}(q^{2k-2j+2};q^2)_{j-1}.
\end{align*}
\end{theorem}
\noindent
Subsequently, the authors carry out the same treatment for the same class of partitions but with the distinct parts being incongruent modulo $2$ with the smallest part. Special cases of these results lead to new identities and inequalities for partitions into distinct parts as stated below.

\begin{theorem}{\cite[Corollary 5(a)]{An-Ba25}}
\mylabel{thm:An-Ba3}
For any integer $n>2$, we have
$$\mathrm{spt}2_d(n)=2p_d(n-1)-p_d(n),$$
where $p_d(n)$ is the number of partitions of $n$ into distinct parts.
\end{theorem}

In this paper, we extend the work in \cite{An-Ba25} to overpartitions. Recall that an overpartition is a partition in which the first occurrence of a number may be overlined \cite{Co-Lo04}. Overpartitions extend classical partitions very naturally and have become a natural framework connecting combinatorics, $q$-series, and modular forms. The foundational work laid by Corteel and Lovejoy \cite{Co-Lo04} developed their combinatorial structures and links to basic hypergeometric series, establishing overpartitions as a central object in the modern theory of partitions. Subsequent results such as overpartition analogs of Rogers-Ramanujan-type identities \cite{Ba-Ye22, Lo04}, the Dyson rank \cite{Lo05, Lo08} and its connection to Maass forms by Bringmann and Lovejoy \cite{Br-Lo07}, among others, show their deep arithmetic and analytic significance.
\\

To start with, we consider a very natural extension of Andrews' and Bachraoui's aforementioned $\mathrm{spt}k_d(n)$ function to overpartitions. 

\begin{definition}
\mylabel{def:sptk}
Let $\overline{\mathrm{Spt}}k(n)$ denote the set of overpartitions of $n$ where the smallest non-overlined part, say $s(\pi)$, appears $k$ times and every overlined part is bigger than $s(\pi)$. Accordingly, let $\overline{\mathrm{spt}}k(n)$ be the cardinality of $\overline{\mathrm{Spt}}k(n)$. 
\\\\
Also, let $a_e(k,n)$ (resp., $a_o(k,n)$) denote the number of partitions $\pi$ of $n$ in $\overline{\mathrm{Spt}}k(n)$ where the number of parts that are greater than $s(\pi)$ is even (resp., odd) and let
\[\overline{\mathrm{spt}}k'(n):=a_e(k,n)-a_o(k,n). \numberthis \label{def;defspt'}\]
\end{definition}
\noindent
Then, by standard combinatorial arguments, we may write
\[\sum\limits_{n=1}^{\infty}\overline{\mathrm{spt}}k(n)q^n=\sum\limits_{n=1}^{\infty}q^{kn}\,\frac{(-q^{n+1};q)_{\infty}}{(q^{n+1};q)_{\infty}}, \numberthis\label{def;spt}\]
and 
\[\sum\limits_{n=1}^{\infty}\overline{\mathrm{spt}}k'(n)q^n=\sum\limits_{n=1}^{\infty}q^{kn}\,\frac{(q^{n+1};q)_{\infty}}{(-q^{n+1};q)_{\infty}}. \numberthis \label{def;spt'}\]

\begin{definition}
\mylabel{def:sptko}
Next, we define $\overline{\mathrm{Spt}}k_o(n)$ to be the set of overpartitions $\pi$ of $n$ where the smallest non-overlined part, say $s(\pi)$, appears $k$ times, every overlined part is larger than $s(\pi)$ and all parts other than $s(\pi)$ are incongruent modulo $2$ with $s(\pi)$. Accordingly, let $\overline{\mathrm{spt}}k_o(n)$ be the cardinality of $\overline{\mathrm{Spt}}k_o(n)$. 
\\\\
Also, let $B_e(k,n)$ (resp., $B_o(k,n)$) denote the set of overpartitions $\pi$ of $n$ in $\overline{\mathrm{Spt}}k_o(n)$ where the number of parts that are greater than $s(\pi)$ is even (resp., odd). Accordingly, let $b_e(k,n)$ (resp., $b_o(k,n)$) be the cardinality of $B_e(k,n)$ (resp., $B_o(k,n)$) and let
\begin{equation}
\label{eq:sptdiff}
\overline{\mathrm{spt}}k_o'(n) := b_e(k,n)-b_o(k,n).
\end{equation} 
\end{definition}
\noindent
Once again, by standard combinatorial arguments, we note that
\[\sum\limits_{n=1}^{\infty}\overline{\mathrm{spt}}k_o(n)q^n=\sum\limits_{n=1}^{\infty}q^{kn}\,\frac{(-q^{n+1};q^2)_{\infty}}{(q^{n+1};q^2)_{\infty}}, \numberthis \label{def;spto}\]
and
\[\sum\limits_{n=1}^{\infty}\overline{\mathrm{spt}}k_o'(n)q^n=\sum\limits_{n=1}^{\infty}q^{kn}\,\frac{(q^{n+1};q^2)_{\infty}}{(-q^{n+1};q^2)_{\infty}}. \numberthis \label{def;spto'}\]
\\
\noindent
The paper is organized as follows. In Section \sect{secresults}, we state our main results concerning the generating functions for the classes of overpartition described previously in the introduction. Section \sect{secproofs} contains proofs of the results stated in Section \sect{secresults}. In Section \ref{sec:secappl}, we present some corollaries which are of interest in their own right. Finally, in Section \ref{sec:secconclusions}, we discuss some open questions arising from our work.

\section{Statement of main results}
\label{sec:secresults}
Our first result is a direct analog of Theorem \thm{An-Ba1} \cite[Theorem 1]{An-Ba25}.

\begin{theorem}
\mylabel{thm:sptk}
For any positive integer $k$, let $\overline{\mathrm{spt}}k(n)$ denote the partition function as in Definition \ref{def:sptk}. Then,
\begin{align*}
\sum\limits_{n=1}^{\infty}\overline{\mathrm{spt}}k(n)q^n=\overline{P}_k(q)\frac{(-q^2;q)_{\infty}}{(q^2;q)_{\infty}}+(-1)^k\frac{(q;q)_{k-1}}{(-q;q)_{k}},
\end{align*}
where 
\begin{align*}
\overline{P}_k(q)=\begin{cases}
1 & \text{if } k=1, \\
\frac{(q^{k-1}-1)\overline{P}_{k-1}(q)+q^{k-1}(1+q)}{1+q^k}& \text{if } k>1.
\end{cases}
\end{align*}
\end{theorem}
\noindent
Next, we give a closed formula for $\overline{P}_k(q)$ as rational functions of $q$, analogous to Theorem \thm{An-Ba2}, \cite[Theorem 2]{An-Ba25}.

\begin{theorem}
\mylabel{thm:sptkpoly} For any positive integer $k$, the rational functions $\overline{P}_{k}(q)$ in the above theorem are given by
\begin{align*}
\overline{P}_{k}(q)=\frac{q^{k-1}(1+q)}{1+q^{k}}-\frac{(q;q)_{k-1}}{(-q^2;q)_{k-1}}\sum\limits_{j=0}^{k-2}(-1)^jq^{k-j-2}\frac{(-q;q)_{k-j-2}}{(q;q)_{k-j-2}},
\end{align*}
where the empty sum in the case of $k=1$ is assumed to be $0$.
\end{theorem}

\begin{theorem}
\mylabel{thm:sptk'}
For any positive integer $k$, the partition function $\overline{\mathrm{spt}}k'(n)$ defined in \eqref{def;defspt'} satisfies
\begin{align*}
\sum\limits_{n=1}^{\infty}\overline{\mathrm{spt}}k'(n)q^n=\frac{(q;q)_{k-1}}{(-q;q)_{k}}-\frac{(q;q)_{\infty}}{(-q;q)_{\infty}}.
\end{align*}
\end{theorem}
\noindent
Our next result is an analog of \cite[Theorem 5]{An-Ba25}.

\begin{theorem}
\mylabel{thm:sptko}
For any positive integer $k$, let $\overline{\mathrm{spt}}k_o(n)$ be as in Definition \ref{def:sptko}. We have
\begin{align*}
\sum\limits_{n=1}^{\infty}\overline{\mathrm{spt}}k_o(n)q^n=\overline{V}_{k}(q)\frac{(-q^4;q^2)_{\infty}}{(q^2;q^2)_{\infty}}+\overline{W}_{k}(q)\frac{(-q;q^2)_{\infty}}{(q^3;q^2)_{\infty}}+2(-q)^k\frac{(q^2;q^2)_{k-1}}{(-q^2;q^2)_{k}},
\end{align*}
where 
\begin{align*}
\overline{V}_{k}(q)=\begin{cases}
2q & \text{if } k=1, \\
\frac{(q^{2k-1}-q)\overline{V}_{k-1}(q)+2q^k(1+q^2)}{1+q^{2k}}& \text{if } k>1,
\end{cases}
\end{align*}
and
\begin{align*}
\overline{W}_{k}(q)=\begin{cases}
\frac{q}{1+q^2} & \text{if } k=1, \\
\frac{(q^{2k-1}-q)\overline{W}_{k-1}(q)+q^{2k-1}}{1+q^{2k}}& \text{if } k>1.
\end{cases}
\end{align*}
\end{theorem}
\noindent
Now, we give closed formulas for $\overline{V}_{k}(q)$ and $\overline{W}_{k}(q)$ as rational functions of $q$, in analogy to \cite[Theorem 6]{An-Ba25}.

\begin{theorem} For a positive integer $k$,
\mylabel{thm:sptkopoly}
\begin{align*} 
\overline{V}_{k}(q)&=(1+q^2)\left(\frac{2q^{k}}{1+q^{2k}}+2q^k\frac{(q^2;q^2)_{k-1}}{(-q^2;q^2)_{k}}\sum\limits_{j=0}^{k-2}(-1)^{j+1}\frac{(-q^2;q^2)_{k-j-2}}{(q^2;q^2)_{k-j-2}}\right),
\\
\overline{W}_{k}(q)&=\frac{q^{2k-1}}{1+q^{2k}}+q^{2k-2}\frac{(q^2;q^2)_{k-1}}{(-q^2;q^2)_{k}}\sum\limits_{j=0}^{k-2}\frac{(-1)^{j+1}}{q^j}\frac{(-q^2;q^2)_{k-j-2}}{(q^2;q^2)_{k-j-2}},
\end{align*}
where the empty sum in the case of $k=1$ is assumed to be $0$.
\end{theorem}
\noindent
Our last main theorem is an analog of \cite[Theorem 3]{An-Ba25}.

\begin{theorem}
\mylabel{thm:sptko'}
For any positive integer $k$, we have
\begin{align*}
\sum\limits_{n=1}^{\infty}\overline{\mathrm{spt}}k_o'(n)q^n=\overline{T}_{k}(q)\frac{(q;q^2)_{\infty}}{(-q;q^2)_{\infty}}+2q^k\frac{(q^2;q^2)_{k-1}}{(-q^2;q^2)_{k}},
\end{align*}
where 
\begin{align*}
\overline{T}_{k}(q)=\begin{cases}
-\frac{q(1+q)}{1+q^2} & \text{if } k=1, \\
\frac{(q-q^{2k-1})\overline{T}_{k-1}(q)-q^{2k-1}(1+q)}{1+q^{2k}} & \text{if } k>1.
\end{cases}
\end{align*}
\end{theorem}
\noindent
Now, we give a closed formula for $\overline{T}_{k}(q)$ as rational functions of $q$ analogous to \cite[Theorem 4]{An-Ba25}.

\begin{theorem} For a positive integer $k$,
\mylabel{thm:sptko'poly}
\begin{align*} 
\overline{T}_{k}(q)=-\frac{q^{2k-1}(1+q)}{1+q^{2k}}-(1+q)\frac{(q^2;q^2)_{k-1}}{(-q^2;q^2)_{k}}\sum\limits_{j=0}^{k-2}q^{2k-j-2}\frac{(-q^2;q^2)_{k-j-2}}{(q^2;q^2)_{k-j-2}},
\end{align*}
where the empty sum in the case of $k=1$ is assumed to be $0$.
\end{theorem}

\section{Proofs}
\label{sec:secproofs}

We begin by recording the $q$-binomial theorem \cite{An98,Ga-Rah04}
\begin{equation}
\mylabel{eq:qbinom}
\sum_{n\geq 0} \frac{(a;q)_n}{(q;q)_n} z^n
=\frac{(az;q)_{\infty}}{(z;q)_{\infty}},
\end{equation}
and the following formula of Andrews, Subbarao, and Vidyasagar \cite{An-Su-Vi72}
\begin{equation}
\mylabel{eq:An-Su-Vi}
\sum_{n\geq 0} \frac{(a;q)_n}{(b;q)_n} q^n
=\frac{q(a;q)_\infty}{b(b;q)_\infty \big(1-\frac{aq}{b}\big)} + \frac{1-\frac{q}{b}}{1-\frac{aq}{b}}
\end{equation}
which we use several times in this section.
\subsection{Proof of Theorem \thm{sptk}}

\begin{proof}
We prove this using the principle of mathematical induction on $k$.
\\\\
The theorem is true for the case $k=1$ as argued below. By the generating function given in \eqref{def;spt}, and letting $a=q$, and $b=-q$ in equation \eqn{An-Su-Vi}, we have
\begin{align*}
\sum\limits_{n=1}^{\infty}\overline{\mathrm{spt}}1(n)q^n&=\sum\limits_{n=1}^{\infty}q^n\,\frac{(-q^{n+1};q)_{\infty}}{(q^{n+1};q)_{\infty}}
\\
&=\frac{(-q;q)_{\infty}}{(q;q)_{\infty}}\sum\limits_{n=1}^{\infty}\frac{(q;q)_{n}}{(-q;q)_{n}}q^n
\\
&=\frac{(-q;q)_{\infty}}{(q;q)_{\infty}}\left(-\frac{(q;q)_{\infty}}{(-q;q)_{\infty}(1+q)}+\frac{2}{1+q}-1\right)
\\
&=\frac{-1}{1+q}+\frac{2\,(-q;q)_{\infty}}{(1+q)(q;q)_{\infty}}-\frac{(-q;q)_{\infty}}{(q;q)_{\infty}}
\\
&=\frac{(-q^2;q)_{\infty}}{(q^2;q)_{\infty}}-\frac{1}{1+q}.
\end{align*}
\\
Next, let us assume that the theorem holds for $k-1$, $k\ge 2$, and we prove the result for $k$. Setting $\overline{C}_k(q):=\sum\limits_{n=1}^{\infty}\overline{\mathrm{spt}}k(n)q^n$, we have
\begin{align*}
\overline{C}_k(q)
&=\sum\limits_{n=1}^{\infty}q^{kn}\,\frac{(-q^{n+1};q)_{\infty}}{(q^{n+1};q)_{\infty}}
\\
&=\sum\limits_{n=1}^{\infty}q^{(k-1)n}\,\left(\frac{1+q^n-1}{1-q^n}\right)\,(1-q^n)\,\frac{(-q^{n+1};q)_{\infty}}{(q^{n+1};q)_{\infty}}
\\
&=-\sum\limits_{n=1}^{\infty}q^{(k-1)n}\,\left(\frac{1-q^n}{1-q^n}\right)\,\frac{(-q^{n+1};q)_{\infty}}{(q^{n+1};q)_{\infty}}
\\
&\hspace{5mm}+\sum\limits_{n=1}^{\infty}q^{(k-1)n}\,\left(\frac{1+q^n}{1-q^n}\right)\,(1-q^n)\,\frac{(-q^{n+1};q)_{\infty}}{(q^{n+1};q)_{\infty}}
\\
&=-\sum\limits_{n=1}^{\infty}q^{(k-1)n}\,\frac{(-q^{n+1};q)_{\infty}}{(q^{n+1};q)_{\infty}}+\sum\limits_{n=1}^{\infty}q^{(k-1)n}\,\frac{(-q^{n};q)_{\infty}}{(q^{n};q)_{\infty}}-\sum\limits_{n=1}^{\infty}q^{kn}\,\frac{(-q^{n};q)_{\infty}}{(q^{n};q)_{\infty}}
\\
&=-\sum\limits_{n=1}^{\infty}q^{(k-1)n}\,\frac{(-q^{n+1};q)_{\infty}}{(q^{n+1};q)_{\infty}}+q^{k-1}\frac{(-q;q)_{\infty}}{(q;q)_{\infty}}+\sum\limits_{n=2}^{\infty}q^{(k-1)n}\,\frac{(-q^{n};q)_{\infty}}{(q^{n};q)_{\infty}}
\\
&\hspace{4mm}-q^{k}\frac{(-q;q)_{\infty}}{(q;q)_{\infty}}-\sum\limits_{n=2}^{\infty}q^{kn}\,\frac{(-q^{n};q)_{\infty}}{(q^{n};q)_{\infty}}
\\
&=-\sum\limits_{n=1}^{\infty}q^{(k-1)n}\,\frac{(-q^{n+1};q)_{\infty}}{(q^{n+1};q)_{\infty}}+q^{k-1}\frac{(-q;q)_{\infty}}{(q;q)_{\infty}}+q^{k-1}\sum\limits_{n=1}^{\infty}q^{(k-1)n}\,\frac{(-q^{n+1};q)_{\infty}}{(q^{n+1};q)_{\infty}}
\\
&\hspace{4mm}-q^{k}\frac{(-q;q)_{\infty}}{(q;q)_{\infty}}-q^k\sum\limits_{n=1}^{\infty}q^{kn}\,\frac{(-q^{n+1};q)_{\infty}}{(q^{n+1};q)_{\infty}}
\\
&=-\left(\overline{P}_{k-1}(q)\frac{(-q^2;q)_{\infty}}{(q^2;q)_{\infty}}+(-1)^{k-1}\frac{(q;q)_{k-2}}{(-q;q)_{k-1}}\right)+q^{k-1}\frac{(-q;q)_{\infty}}{(q;q)_{\infty}}
\\
&\hspace{5mm}+q^{k-1}\left(\overline{P}_{k-1}(q)\frac{(-q^2;q)_{\infty}}{(q^2;q)_{\infty}}+(-1)^{k-1}\frac{(q;q)_{k-2}}{(-q;q)_{k-1}}\right)-q^{k}\frac{(-q;q)_{\infty}}{(q;q)_{\infty}}
\\
&\hspace{5mm}-q^k\overline{C}_{k}(q)
\\
&=(q^{k-1}-1)\overline{P}_{k-1}(q)\frac{(-q^2;q)_{\infty}}{(q^2;q)_{\infty}}+(q^{k-1}-q^k)\frac{(-q;q)_{\infty}}{(q;q)_{\infty}}
\\
&\hspace{5mm}+(-1)^k(1-q^{k-1})\frac{(q;q)_{k-2}}{(-q;q)_{k-1}}-q^k\overline{C}_{k}(q)
\\
&=\left((q^{k-1}-1)\overline{P}_{k-1}(q)+\frac{(q^{k-1}-q^k)(1+q)}{1-q}\right)\frac{(-q^2;q)_{\infty}}{(q^2;q)_{\infty}}
\\
&\hspace{5mm}+(-1)^k\frac{(q;q)_{k-1}}{(-q;q)_{k-1}}-q^k\overline{C}_k(q).
\end{align*}
Adding $q^k\overline{C}_k(q)$ on both sides gives us the theorem.
\end{proof}

\subsection{Proof of Theorem \thm{sptkpoly}}

\begin{proof}
Let $(\overline{p}_k(q))=(\overline{p}_k(q))_{k\geq 1}$ be the sequence defined by 
$$\overline{p}_k(q)=\frac{q^{k-1}(1+q)}{1+q^{k}}-\frac{(q;q)_{k-1}}{(-q^2;q)_{k-1}}\sum\limits_{j=0}^{k-2}(-1)^jq^{k-j-2}\frac{(-q;q)_{k-j-2}}{(q;q)_{k-j-2}}.$$
\\
First we note that $$\overline{p}_1(q)=1,$$
which is the same as $\overline{P}_1(q)$ in Theorem \thm{sptk}.
\\\\
Next we show below that $\overline{p}_k(q)$ satisfies the same recurrence as $\overline{P}_k(q)$ given in Theorem \thm{sptk}.
\begin{align*}
\overline{p}_k(q)&=\frac{q^{k-1}(1+q)}{1+q^k}-\frac{(q;q)_{k-1}}{(-q^2;q)_{k-1}}\left(q^{k-2}\frac{(-q;q)_{k-2}}{(q;q)_{k-2}}+\sum_{j=1}^{k-2}(-1)^j q^{k-j-2}\frac{(-q;q)_{k-j-2}}{(q;q)_{k-j-2}}\right)
\\
&= \frac{q^{k-1}(1+q)}{1+q^k}-\frac{q^{k-2}(1-q^{k-1})(1+q)}{(1+q^{k-1})(1+q^k)}+\frac{(q;q)_{k-1}}{(-q^2;q)_{k-1}}\sum_{m=0}^{k-3}(-1)^m q^{k-m-3}\frac{(-q;q)_{k-m-3}}{(q;q)_{k-m-3}}
\\
&= \frac{q^{k-1}(1+q)}{1+q^k}-\frac{q^{k-2}(1-q^{k-1})(1+q)}{(1+q^{k-1})(1+q^k)}
\\
&\hspace{5mm}-\frac{-1+q^{k-1}}{1+q^k}\frac{(q;q)_{k-2}}{(-q^2;q)_{k-2}}\sum_{m=0}^{k-3}(-1)^m q^{k-m-3}
\frac{(-q;q)_{k-m-3}}{(q;q)_{k-m-3}}
\\
&= \frac{q^{k-1}(1+q)}{1+q^k}-\frac{q^{k-2}(1-q^{k-1})(1+q)}{(1+q^{k-1})(1+q^k)}+\frac{-1+q^{k-1}}{1+q^k}
\left(\overline{p}_{k-1}(q)-\frac{q^{k-2}(1+q)}{1+q^{k-1}}\right)
\\
&=\frac{q^{k-1}(1+q)+(q^{k-1}-1)\overline{p}_{k-1}(q)}{1+q^k}.
\end{align*}
This completes the proof of the theorem.
\end{proof}

\subsection{Proof of Theorem \thm{sptk'}}

\begin{proof}
By the generating function given in \eqref{def;spt'}, we have 
\begin{align*}
\sum\limits_{n=1}^{\infty}\overline{\mathrm{spt}}k'(n)q^n&=\sum\limits_{n=1}^{\infty}q^{kn}\,\frac{(q^{n+1};q)_{\infty}}{(-q^{n+1};q)_{\infty}}
\\
&=\frac{(q;q)_{\infty}}{(-q;q)_{\infty}}\sum\limits_{n=1}^{\infty}\frac{(-q;q)_{n}}{(q;q)_{n}}q^{kn}
\\
&=\frac{(q;q)_{\infty}}{(-q;q)_{\infty}}\frac{(-q^{k+1};q)_{\infty}}{(q^k;q)_{\infty}}-\frac{(q;q)_{\infty}}{(-q;q)_{\infty}}
\\
&=\frac{(q;q)_{k-1}}{(-q;q)_{k}}-\frac{(q;q)_{\infty}}{(-q;q)_{\infty}},\\
\end{align*}
where in the third step we have used the $q$-binomial theorem in Equation \eqn{qbinom} with $a=-q$, and $z=q^k$.
\end{proof}

\subsection{Proof of Theorem \thm{sptko}}

\begin{proof}
We again prove this using mathematical induction on $k$.
The theorem is true for the case $k=1$ as argued below. The generating function for $\overline{\mathrm{spt}}k_o(n)$ is given in \eqref{def;spto}. Letting $q \rightarrow q^2$, $a=q$ and $b=-q$ in Equation \eqn{An-Su-Vi} for the first sum below and $q \rightarrow q^2$, $a=q^2$ and $b=-q^2$ in Equation \eqn{An-Su-Vi} for the second sum below, we have
\begin{align*}
\sum\limits_{n=1}^{\infty}\overline{\mathrm{spt}}1_o(n)q^n&=\sum\limits_{n=1}^{\infty}q^n\,\frac{(-q^{n+1};q^2)_{\infty}}{(q^{n+1};q^2)_{\infty}}
\\
&=\sum\limits_{n=1}^{\infty}q^{2n}\,\frac{(-q^{2n+1};q^2)_{\infty}}{(q^{2n+1};q^2)_{\infty}}+\sum\limits_{n=0}^{\infty}q^{2n+1}\,\frac{(-q^{2n+2};q^2)_{\infty}}{(q^{2n+2};q^2)_{\infty}}
\\
&=\frac{(-q;q^2)_{\infty}}{(q;q^2)_{\infty}}\sum\limits_{n=1}^{\infty}\frac{(q;q^2)_{n}}{(-q;q^2)_{n}}q^{2n}+q\frac{(-q^2;q^2)_{\infty}}{(q^2;q^2)_{\infty}}\sum\limits_{n=0}^{\infty}\frac{(q^2;q^2)_{n}}{(-q^2;q^2)_{n}}q^{2n}
\\
&=\frac{(-q;q^2)_{\infty}}{(q;q^2)_{\infty}}\left(-\frac{q(q;q^2)_{\infty}}{(-q;q^2)_{\infty}(1+q^2)}+\frac{1+q}{1+q^2}-1\right)
\\
&\hspace{5mm}+\frac{q(-q^2;q^2)_{\infty}}{(q^2;q^2)_{\infty}}\left(-\frac{(q^2;q^2)_{\infty}}{(-q^2;q^2)_{\infty}(1+q^2)}+\frac{2}{1+q^2}\right)
\\
&=\frac{-2q}{1+q^2}+\frac{2q}{1+q^2}\frac{(-q^2;q^2)_{\infty}}{(q^2;q^2)_{\infty}}+\frac{q(1-q)}{1+q^2}\frac{(-q;q^2)_{\infty}}{(q;q^2)_{\infty}}
\\
&=2q\frac{(-q^4;q^2)_{\infty}}{(q^2;q^2)_{\infty}}+\frac{q}{1+q^2}\frac{(-q;q^2)_{\infty}}{(q^3;q^2)_{\infty}}+\frac{-2q}{1+q^2}.
\end{align*}
\\
Next, let us assume that the theorem holds for $k-1$. Then setting $\overline{U}_k(q):=\sum\limits_{n=1}^{\infty}\overline{\mathrm{spt}}k_o(n)q^n$, we have
\begin{align*}
\overline{U}_k(q)&=\sum\limits_{n=1}^{\infty}q^{kn}\,\frac{(-q^{n+1};q^2)_{\infty}}{(q^{n+1};q^2)_{\infty}}
\\
&=q^{k}\,\frac{(-q^{2};q^2)_{\infty}}{(q^{2};q^2)_{\infty}}+\sum\limits_{n=2}^{\infty}q^{kn}\,\frac{(-q^{n+1};q^2)_{\infty}}{(q^{n+1};q^2)_{\infty}}
\\
&=q^{k}\,\frac{(-q^{2};q^2)_{\infty}}{(q^{2};q^2)_{\infty}}+\sum\limits_{n=2}^{\infty}q^{(k-1)n}\,\frac{-q+q(1+q^{n-1})}{1-q^{n-1}}\,(1-q^{n-1})\,\frac{(-q^{n+1};q^2)_{\infty}}{(q^{n+1};q^2)_{\infty}}
\\
&=q^{k}\,\frac{(-q^{2};q^2)_{\infty}}{(q^{2};q^2)_{\infty}}-q\sum\limits_{n=2}^{\infty}q^{(k-1)n}\,\frac{(-q^{n+1};q^2)_{\infty}}{(q^{n+1};q^2)_{\infty}}+q\sum\limits_{n=2}^{\infty}q^{(k-1)n}\,\frac{(-q^{n-1};q^2)_{\infty}}{(q^{n-1};q^2)_{\infty}}
\\
&-\sum\limits_{n=2}^{\infty}q^{kn}\,\frac{(-q^{n-1};q^2)_{\infty}}{(q^{n-1};q^2)_{\infty}}
\\
&=2q^{k}\,\frac{(-q^{2};q^2)_{\infty}}{(q^{2};q^2)_{\infty}}-q\overline{U}_{k-1}(q)+q^{2k-1}\overline{U}_{k-1}(q)+q^{2k-1}\frac{(-q;q^2)_{\infty}}{(q;q^2)_{\infty}}-q^{2k}\overline{U}_{k}(q)
\\
&-q^{2k}\frac{(-q;q^2)_{\infty}}{(q;q^2)_{\infty}}
\\
&=2q^k(1+q^2)\frac{(-q^4;q^2)_{\infty}}{(q^2;q^2)_{\infty}}+(q^{2k-1}-q)\Bigg(\overline{V}_{k-1}(q)\frac{(-q^4;q^2)_{\infty}}{(q^2;q^2)_{\infty}}+\overline{W}_{k-1}(q)\frac{(-q;q^2)_{\infty}}{(q^3;q^2)_{\infty}}
\\
&+2(-q)^{k-1}\frac{(q^2;q^2)_{k-2}}{(-q^2;q^2)_{k-1}}\Bigg)+(q^{2k-1}-q^{2k})\frac{(-q;q^2)_{\infty}}{(q;q^2)_{\infty}}-q^{2k}\overline{U}_k(q)
\\
&=\frac{(-q^4;q^2)_{\infty}}{(q^2;q^2)_{\infty}}\Big(2q^k(1+q^2)+(q^{2k-1}-q)\overline{V}_{k-1}(q)\Big)
\\
&+\frac{(-q;q^2)_{\infty}}{(q^3;q^2)_{\infty}}\Big(q^{2k-1}+(q^{2k-1}-q)\overline{W}_{k-1}(q)\Big)
\\
&+2(q^{2k-1}-q)(-q)^{k-1}\frac{(q^2;q^2)_{k-2}}{(-q^2;q^2)_{k-1}}-q^{2k}\overline{U}_k(q).
\end{align*}
Adding $q^{2k}\overline{U}_k(q)$ on both sides gives us our theorem.
\end{proof}

\subsection{Proof of Theorem \thm{sptkopoly}}

\begin{proof}
Let $(\overline{v}_k(q))=(\overline{v}_k(q))_{k\geq 1}$ be the sequence defined by 
$$\overline{v}_k(q)=(1+q^2)\left(\frac{2q^{k}}{1+q^{2k}}+2q^k\frac{(q^2;q^2)_{k-1}}{(-q^2;q^2)_{k}}\sum\limits_{j=0}^{k-2}(-1)^{j+1}\frac{(-q^2;q^2)_{k-j-2}}{(q^2;q^2)_{k-j-2}}\right).$$
\\
First we note that $$\overline{v}_1(q)=2q,$$
which is the same as $\overline{V}_1(q)$ in Theorem \thm{sptko}.
\\\\
Next we show below that $\overline{v}_k(q)$ satisfies the same recurrence as $\overline{V}_k(q)$ given in Theorem \thm{sptko}.
\begin{align*}
\overline{v}_k(q)&=\frac{2q^k(1+q^2)}{1+q^{2k}}+\frac{2q^k(1+q^2)(q^2;q^2)_{k-1}}{(-q^2;q^2)_k}\left(\frac{-(-q^2;q^2)_{k-2}}{(q^2;q^2)_{k-2}}+ \sum_{j=1}^{k-2}(-1)^{j+1}\frac{(-q^2;q^2)_{k-j-2}}{(q^2;q^2)_{k-j-2}}\right)
\\
&=\frac{2q^k(1+q^2)}{1+q^{2k}}-\frac{2q^k(1+q^2)(1-q^{2k-2})}{(1+q^{2k})(1+q^{2k-2})}
\\
&\hspace{5mm}+\frac{2q^k(1+q^2)(q^2;q^2)_{k-1}}{(-q^2;q^2)_k}\sum_{j=1}^{k-2}(-1)^{j+1}\frac{(-q^2;q^2)_{k-j-2}}{(q^2;q^2)_{k-j-2}}
\\
&=\frac{2q^k(1+q^2)}{1+q^{2k}}-\frac{2q^k(1+q^2)(1-q^{2k-2})}{(1+q^{2k})(1+q^{2k-2})}
\\
&\hspace{5mm}+\frac{2q^{k-1}(1+q^2)(q^2;q^2)_{k-2}}{(-q^2;q^2)_{k-1}}q\,\frac{1-q^{2k-2}}{1+q^{2k}}\sum_{m=0}^{k-3}(-1)^m\frac{(-q^2;q^2)_{(k-1)-m-2}}{(q^2;q^2)_{(k-1)-m-2}}
\\
&=\frac{2q^k(1+q^2)}{1+q^{2k}}-\frac{2q^k(1+q^2)(1-q^{2k-2})}{(1+q^{2k})(1+q^{2k-2})}
\\
&\hspace{5mm}-q\,\frac{1-q^{2k-2}}{1+q^{2k}}\left(\frac{2q^{k-1}(1+q^2)(q^2;q^2)_{k-2}}{(-q^2;q^2)_{k-1}}\sum_{m=0}^{k-3}(-1)^{m+1}\frac{(-q^2;q^2)_{(k-1)-m-2}}{(q^2;q^2)_{(k-1)-m-2}}\right)
\\
&= \frac{2q^k(1+q^2)}{1+q^{2k}}-\frac{2q^k(1+q^2)(1-q^{2k-2})}{(1+q^{2k})(1+q^{2k-2})}-\frac{q(1-q^{2k-2})}{1+q^{2k}}\left(\overline{v}_{k-1}(q)-\frac{2q^{k-1}(1+q^2)}{1+q^{2k-2}}\right)
\\
&= \frac{2q^k(1+q^2)}{1+q^{2k}}-\frac{2q^k(1+q^2)(1-q^{2k-2})}{(1+q^{2k})(1+q^{2k-2})}+\frac{q^{2k-1}-q}{1+q^{2k}}\,\overline{v}_{k-1}(q)
\\
&\hspace{5mm}+\frac{2q^k(1-q^{2k-2})(1+q^2)}{(1+q^{2k})(1+q^{2k-2})}
\\
&=\frac{(q^{2k-1}-q)\overline{v}_{k-1}(q)+2q^k(1+q^2)}{1+q^{2k}}.
\end{align*}
\\
The same can be done to prove the formula for $\overline{W}_{k}(q)$ and the steps of the proof follow very closely. We leave it as an exercise to the reader. This proves the theorem.
\end{proof}

\subsection{Proof of Theorem \thm{sptko'}}

\begin{proof}
Once again, we prove this using mathematical induction on $k$.
\\\\
The theorem is true for the case $k=1$ as argued below. The generating function for $\overline{\mathrm{spt}}k_o'(n)$ is given in \eqref{def;spto'}. Letting $q \rightarrow q^2$, $a=-q$ and $b=q$ in Equation \eqn{An-Su-Vi} for the first sum below and $q \rightarrow q^2$, $a=-q^2$ and $b=q^2$ in Equation \eqn{An-Su-Vi} for the second sum below, we have
\begin{align*}
\sum\limits_{n=1}^{\infty}\overline{\mathrm{spt}}1_o'(n)q^n&=\sum\limits_{n=1}^{\infty}q^n\,\frac{(q^{n+1};q^2)_{\infty}}{(-q^{n+1};q^2)_{\infty}}
\\
&=\sum\limits_{n=1}^{\infty}q^{2n}\,\frac{(q^{2n+1};q^2)_{\infty}}{(-q^{2n+1};q^2)_{\infty}}+\sum\limits_{n=0}^{\infty}q^{2n+1}\,\frac{(q^{2n+2};q^2)_{\infty}}{(-q^{2n+2};q^2)_{\infty}}
\\
&=\frac{(q;q^2)_{\infty}}{(-q;q^2)_{\infty}}\sum\limits_{n=1}^{\infty}\frac{(-q;q^2)_{n}}{(q;q^2)_{n}}q^{2n}+q\frac{(q^2;q^2)_{\infty}}{(-q^2;q^2)_{\infty}}\sum\limits_{n=0}^{\infty}\frac{(-q^2;q^2)_{n}}{(q^2;q^2)_{n}}q^{2n}
\\
&=\frac{(q;q^2)_{\infty}}{(-q;q^2)_{\infty}}\left(\frac{q(-q;q^2)_{\infty}}{(q;q^2)_{\infty}(1+q^2)}+\frac{1-q}{1+q^2}-1\right)
\\
&\hspace{5mm}+\frac{q(q^2;q^2)_{\infty}}{(-q^2;q^2)_{\infty}}\left(\frac{(-q^2;q^2)_{\infty}}{(q^2;q^2)_{\infty}(1+q^2)}+\frac{0}{1+q^2}\right)
\\
&=\frac{2q}{1+q^2}+\frac{-q(1+q)}{1+q^2}\frac{(q;q^2)_{\infty}}{(-q;q^2)_{\infty}}.
\end{align*}
\\
Next, let us assume that the theorem holds for $k-1$. Then setting $\overline{A}_k(q):=\sum\limits_{n=1}^{\infty}\overline{\mathrm{spt}}k_o'(n)q^n$, we have
\begin{align*}
\overline{A}_k(q)&=\sum\limits_{n=1}^{\infty}q^{kn}\,\frac{(q^{n+1};q^2)_{\infty}}{(-q^{n+1};q^2)_{\infty}}
\\
&=\sum\limits_{n=1}^{\infty}q^{(k-1)n}\,\frac{q-q(1-q^{n-1})}{1+q^{n-1}}\,(1+q^{n-1})\,\frac{(q^{n+1};q^2)_{\infty}}{(-q^{n+1};q^2)_{\infty}}
\\
&=q\sum\limits_{n=1}^{\infty}q^{(k-1)n}\,\frac{(q^{n+1};q^2)_{\infty}}{(-q^{n+1};q^2)_{\infty}}-q\sum\limits_{n=1}^{\infty}q^{(k-1)n}\,\frac{(q^{n-1};q^2)_{\infty}}{(-q^{n-1};q^2)_{\infty}}-\sum\limits_{n=1}^{\infty}q^{kn}\,\frac{(q^{n-1};q^2)_{\infty}}{(-q^{n-1};q^2)_{\infty}}
\\
&=q\overline{A}_{k-1}(q)-q\sum\limits_{n=-1}^{\infty}q^{(k-1)(n+2)}\,\frac{(q^{n+1};q^2)_{\infty}}{(-q^{n+1};q^2)_{\infty}}-\sum\limits_{n=-1}^{\infty}q^{k(n+2)}\,\frac{(q^{n+1};q^2)_{\infty}}{(-q^{n+1};q^2)_{\infty}}
\\
&=q\overline{A}_{k-1}(q)-q^{2k-1}\,\frac{(q;q^2)_{\infty}}{(-q;q^2)_{\infty}}-q^{2k}\,\frac{(q;q^2)_{\infty}}{(-q;q^2)_{\infty}}-q^{2k-1}\,\overline{A}_{k-1}(q)-q^{2k}\,\overline{A}_k(q)
\\
&=(q-q^{2k-1})\overline{A}_{k-1}(q)-(q^{2k-1}+q^{2k})\frac{(q;q^2)_{\infty}}{(-q;q^2)_{\infty}}-q^{2k}\overline{A}_k(q)
\\
&=(q-q^{2k-1})\left(\overline{T}_{k-1}(q)\frac{(q;q^2)_{\infty}}{(-q;q^2)_{\infty}}+2q^{k-1}\frac{(q^2;q^2)_{k-2}}{(-q^2;q^2)_{k-1}}\right)-(q^{2k-1}+q^{2k})\frac{(q;q^2)_{\infty}}{(-q;q^2)_{\infty}}
\\
&\hspace{5mm}-q^{2k}\overline{A}_k(q)
\\
&=\left(\overline{T}_{k-1}(q)(q-q^{2k-1})-(q^{2k-1}+q^{2k})\right)\frac{(q;q^2)_{\infty}}{(-q;q^2)_{\infty}}+2q^k\frac{(q^2;q^2)_{k-1}}{(-q^2;q^2)_{k-1}}-q^{2k}\overline{A}_k(q).
\end{align*}
\\
Adding $q^{2k}\overline{A}_k(q)$ on both sides gives us our theorem.
\end{proof}

\subsection{Proof of Theorem \thm{sptko'poly}}

\begin{proof}
Let $(\overline{t}_k(q))=(\overline{t}_k(q))_{k\geq 1}$ be the sequence defined by 
$$\overline{t}_k(q)=-\frac{q^{2k-1}(1+q)}{1+q^{2k}}-(1+q)\frac{(q^2;q^2)_{k-1}}{(-q^2;q^2)_{k}}\sum\limits_{j=0}^{k-2}q^{2k-j-2}\frac{(-q^2;q^2)_{k-j-2}}{(q^2;q^2)_{k-j-2}}.$$
\\
First we note that $$\overline{t}_1(q)=-\frac{q(1+q)}{1+q^2},$$
which is the same as $\overline{T}_1(q)$ in Theorem \thm{sptko'}. Next we show below that $\overline{t}_k(q)$ satisfies the same recurrence as $\overline{T}_k(q)$ given in Theorem \thm{sptko'}.
\begin{align*}
\overline{t}_k(q)&=-\frac{q^{2k-1}(1+q)}{1+q^{2k}}-(1+q)\frac{(q^2;q^2)_{k-1}}{(-q^2;q^2)_{k}}\left(q^{2k-2}
\frac{(-q^2;q^2)_{k-2}}{(q^2;q^2)_{k-2}}+\sum_{j=1}^{k-2}q^{2k-j-2}\frac{(-q^2;q^2)_{k-j-2}}{(q^2;q^2)_{k-j-2}}
\right)
\\
&=-\frac{q^{2k-1}(1+q)}{1+q^{2k}}-\frac{q^{2k-2}(1-q^{2k-2})(1+q)}{(1+q^{2k-2})(1+q^{2k})}
\\
&\hspace{5mm}-(1+q)\frac{(q^2;q^2)_{k-1}}{(-q^2;q^2)_{k}}\sum_{m=0}^{k-3}q^{2k-m-3}\frac{(-q^2;q^2)_{k-m-3}}{(q^2;q^2)_{k-m-3}}
\\
&=-\frac{q^{2k-1}(1+q)}{1+q^{2k}}-\frac{q^{2k-2}(1-q^{2k-2})(1+q)}{(1+q^{2k-2})(1+q^{2k})}
\\
&\hspace{5mm}+q\,\frac{1-q^{2k-2}}{1+q^{2k}}\left(-(1+q)\frac{(q^{2};q^{2})_{k-2}}{(-q^{2};q^{2})_{k-1}}\sum_{m=0}^{k-3}q^{\,2(k-1)-m-2}\frac{(-q^{2};q^{2})_{(k-1)-m-2}}{(q^{2};q^{2})_{(k-1)-m-2}}\right)
\\
&=-\frac{q^{2k-1}(1+q)}{1+q^{2k}}-\frac{q^{2k-2}(1-q^{2k-2})(1+q)}{(1+q^{2k-2})(1+q^{2k})}+q\,\frac{1-q^{2k-2}}{1+q^{2k}}\left(\overline{t}_{k-1}(q)+q^{2k-3}\frac{1+q}{1+q^{2k-2}}\right)
\\
&=\frac{(q-q^{2k-1})\overline{t}_{k-1}(q)-q^{2k-1}(1+q)}{1+q^{2k}}.
\end{align*}
This proves the theorem.
\end{proof}

\section{Applications to overpartitions}
\label{sec:secappl}

The generating function for the number of overpartitions of $n$ denoted by $\overline{p}(n)$ is given by 
$$\sum\limits_{n=1}^{\infty}\overline{p}(n)q^n=\frac{(-q;q)_{\infty}}{(q;q)_{\infty}}.$$ 
\\
And, the generating functions for the number of overpartitions of $n$ into even (resp. odd) parts denoted by $\overline{p}_e(n)$ (resp. $\overline{p}_o(n)$) is respectively given by 
$$\sum\limits_{n=1}^{\infty}\overline{p}_e(n)q^n=\frac{(-q^2;q^2)_{\infty}}{(q^2;q^2)_{\infty}},$$
and 
$$\sum\limits_{n=1}^{\infty}\overline{p}_o(n)q^n=\frac{(-q;q^2)_{\infty}}{(q;q^2)_{\infty}}.$$
\\
Using Theorems \thm{sptk} and \thm{sptko}, one can write $\overline{\mathrm{spt}}k(n)$ and $\overline{\mathrm{spt}}k_o(n)$  as linear combinations of certain subclasses of overpartitions. Below are some examples. 
\\\\
For $k=1$, Theorem \thm{sptk} gives 
$$\sum\limits_{n=1}^{\infty}\overline{\mathrm{spt}}1(n)q^n=\frac{(-q^2;q)_{\infty}}{(q^2;q)_{\infty}}-\frac{1}{1+q}.$$
\\
Rearranging and simplifying leads to 
$$\sum\limits_{n=1}^{\infty}\overline{\mathrm{spt}}1(n)q^n+\sum\limits_{n=1}^{\infty}\overline{\mathrm{spt}}1(n)q^{n+1}+1=\frac{(-q;q)_{\infty}}{(q^2;q)_{\infty}}.$$
\\
Now the following result follows immediately.
\begin{cor}
\label{cor:corsptk}
For $n>1$, $\overline{\mathrm{spt}}1(n)+\overline{\mathrm{spt}}1(n-1)$ equals the number of overpartitions of $n$ with no non-overlined $1$'s.
\end{cor}
\noindent

The cases for higher $k\geq 2$ also lead to similar results but they involve multiple combinations of the spt-functions and corresponding subclasses of overpartitions. Below, we state the one for $\overline{\mathrm{spt}}2(n)$ but do not pursue the same for the other spt-functions. Substituting $k=2$ in Theorem \thm{sptk} and simplifying leads to 
\begin{cor}
\label{cor:corsptk2}
For $n>1$, we have $$\overline{\mathrm{spt}}2(n)+\overline{\mathrm{spt}}2(n-1)+\overline{\mathrm{spt}}2(n-2)+\overline{\mathrm{spt}}2(n-3)+\overline{p}_{uu}(n)=2\overline{p}_{u}(n-1),$$
where $\overline{p}_{uu}(n)$ denotes the number of overpartitions of $n$ with no non-overlined $1$'s and $2$'s and $\overline{p}_{u}(n)$ denotes the number of overpartitions of $n$ with no non-overlined $1$'s.
\end{cor}
\noindent
Next, for $k=1$, Theorem \thm{sptko} gives 
$$\sum\limits_{n=1}^{\infty}\overline{\mathrm{spt}}1_o(n)q^n=2q\frac{(-q^4;q^2)_{\infty}}{(q^2;q^2)_{\infty}}+\frac{q}{1+q^2}\frac{(-q;q^2)_{\infty}}{(q^3;q^2)_{\infty}}+\frac{-2q}{1+q^2}.$$
Rearranging and simplifying leads to 
$$\sum\limits_{n=1}^{\infty}\overline{\mathrm{spt}}1_o(n)q^n+\sum\limits_{n=1}^{\infty}\overline{\mathrm{spt}}1_o(n)q^{n+2}+2q=2q\frac{(-q^2;q^2)_{\infty}}{(q^2;q^2)_{\infty}}+q\frac{(-q;q^2)_{\infty}}{(q^3;q^2)_{\infty}}.$$
This immediately gives the following result.
\begin{cor}
\label{cor:corsptko}
For $n>2$, we have $$\overline{\mathrm{spt}}1_o(n)+\overline{\mathrm{spt}}1_o(n-2)=2\overline{p}_e(n-1)+\overline{p}_{ou}(n-1),$$
where $\overline{p}_{ou}(n)$ denotes the number of overpartitions of $n$ into odd parts with no non-overlined $1$'s.
\end{cor}
\noindent
Lastly, to obtain a similar result for the overpartition difference function discussed in \eqn{sptdiff}, we once again consider the overpartitions of $n$ into odd parts with no non-overlined $1$'s enumerated by $\overline{p}_{ou}(n)$. Then we have 
$$\frac{(q;q^2)_{\infty}}{(-q^3;q^2)_{\infty}}=\sum\limits_{n=1}^{\infty}\overline{p}_{ou}'(n)q^n,$$
where $\overline{p}_{ou}'(n)$ denotes the number of overpartitions counted by $\overline{p}_{ou}(n)$ where the number of parts is even minus the number of overpartitions counted by $\overline{p}_{ou}(n)$ where the number of parts is odd. Now, for $k=1$, Theorem \thm{sptko'} gives 
$$\sum\limits_{n=1}^{\infty}\overline{\mathrm{spt}}1_o'(n)q^n=\frac{2q}{1+q^2}+\frac{-q(1+q)}{1+q^2}\frac{(q;q^2)_{\infty}}{(-q;q^2)_{\infty}}.$$
Rearranging and simplifying leads to 
$$\sum\limits_{n=1}^{\infty}\overline{\mathrm{spt}}1_o'(n)q^n+\sum\limits_{n=1}^{\infty}\overline{\mathrm{spt}}1_o'(n)q^{n+2}=2q-q\frac{(q;q^2)_{\infty}}{(-q^3;q^2)_{\infty}}.$$
\noindent
Below are two examples in tabular form illustrating Corollaries \corol{corsptk} and \corol{corsptko}.
\begin{table}[H]
\small
\centering
\[
\begin{array}{|c|c|c|c|c|}
\hline
\substack{\text{Overpartitions}\\ \text{of }4}
& \in \overline{\mathrm{Spt}}_1(4)
& \substack{\text{Overpartitions}\\ \text{of }3}
& \in \overline{\mathrm{Spt}}_1(3)
& \substack{\text{No non-overlined }1\text{'s}\\ \text{(for }n=4\text{)}} \\
\hline
4 & \checkmark & 3 & \checkmark & \checkmark \\
\overline{4} &  & \overline{3} &  & \checkmark \\
3+1 & \checkmark & 2+1 & \checkmark &  \\
\overline{3}+1 & \checkmark & \overline{2}+1 & \checkmark &  \\
3+\overline{1} &  & 2+\overline{1} &  & \checkmark \\
\overline{3}+\overline{1} &  & \overline{2}+\overline{1} &  & \checkmark \\
2+2 &  & 1+1+1 &  & \checkmark \\
\overline{2}+2 &  & \overline{1}+1+1 &  & \checkmark \\
2+1+1 &  &  &  &  \\
\overline{2}+1+1 &  &  &  &  \\
2+\overline{1}+1 &  &  &  &  \\
1+1+1+1 &  &  &  &  \\
\overline{1}+1+1+1 &  &  &  &  \\
\hline
\end{array}
\]
\caption{Example illustrating Corollary \corol{corsptk} for $n=4$ and $k=1$}
\end{table}

\begin{table}[H]
\small
\centering
\setlength{\arraycolsep}{4pt} 
\[
\begin{array}{|c|c|c|c|c|c|c|}
\hline
\substack{\text{Overpartitions}\\ \text{of }5}
& \in \overline{\mathrm{Spt}}1_o(5)
& \substack{\text{Overpartitions}\\ \text{of }3}
& \in \overline{\mathrm{Spt}}1_o(3)
& \substack{\text{Overpartitions}\\ \text{of }4}
& \substack{\text{All parts even}\\ (n=4)\\ }
& \substack{\text{All parts odd,}\\
\text{no unoverlined}\\
\text{$1$'s}\\
\text{(for }n=4\text{)}}
\\
\hline
5 & \checkmark & 3 & \checkmark & 4 & \checkmark &  \\
\overline{5} &  & \overline{3} &  & \overline{4} & \checkmark &  \\
4+1 & \checkmark & 2+1 & \checkmark & 3+1 &  &  \\
\overline{4}+1 & \checkmark & \overline{2}+1 & \checkmark & \overline{3}+1 &  &  \\
4+\overline{1} &  & 2+\overline{1} &  & 3+\overline{1} &  & \checkmark \\
\overline{4}+\overline{1} &  & \overline{2}+\overline{1} &  & \overline{3}+\overline{1} &  & \checkmark \\
3+2 & \checkmark & 1+1+1 &  & 2+2 & \checkmark &  \\
\overline{3}+2 & \checkmark & \overline{1}+1+1 &  & \overline{2}+2 & \checkmark &  \\
3+\overline{2} &  &  &  & 2+1+1 &  &  \\
\overline{3}+\overline{2} &  &  &  & \overline{2}+1+1 &  &  \\
3+1+1 &  &  &  & 2+\overline{1}+1 &  &  \\
\overline{3}+1+1 &  &  &  & 1+1+1+1 &  &  \\
3+\overline{1}+1 &  &  &  & \overline{1}+1+1+1 &  &  \\
\overline{3}+\overline{1}+1 &  &  &  &  &  &  \\
2+2+1 & \checkmark &  &  &  &  &  \\
\overline{2}+2+1 & \checkmark &  &  &  &  &  \\
2+2+\overline{1} &  &  &  &  &  &  \\
\overline{2}+2+\overline{1} &  &  &  &  &  &  \\
2+1+1+1 &  &  &  &  &  &  \\
\overline{2}+1+1+1 &  &  &  &  &  &  \\
2+\overline{1}+1+1 &  &  &  &  &  &  \\
\overline{2}+\overline{1}+1+1 &  &  &  &  &  &  \\
1+1+1+1+1 &  &  &  &  &  &  \\
\overline{1}+1+1+1+1 &  &  &  &  &  &  \\
\hline
\end{array}
\]
\setlength{\arraycolsep}{6pt} 
\caption{Example illustrating Corollary \corol{corsptko} for $n=5$ and $k=1$}
\end{table}
\noindent
From Table 1, one can see that $\overline{\mathrm{spt}}1(4)+\overline{\mathrm{spt}}1(3)$ = 6 = number of overpartitions of $4$ with no non-overlined $1$'s, which agrees with Corollary \corol{corsptk}.
\\\\
From Table 2, one can see that 
$$\overline{\mathrm{spt}}1_o(5)+\overline{\mathrm{spt}}1_o(3)=10=2\overline{p}_e(4)+\overline{p}_{ou}(4),$$
which agrees with Corollary \corol{corsptko}.
\\
\begin{cor}
\label{cor:corsptko'}
For $n>2$, $$\overline{\mathrm{spt}}1_o'(n)+\overline{\mathrm{spt}}1_o'(n-2)=-\overline{p}_{ou}'(n-1).$$
\end{cor}
\noindent
To illustrate the above corollary, we observe that the overpartitions in $B_e(1,5)$ are $5, 2+2+1, \overline{2}+2+1$ and those in $B_o(1,5)$ are $4+1, \overline{4}+1, 3+2, \overline{3}+2$. Similarly, the sole overpartition in $B_e(1,3)$ is $3$ and the ones in $B_e(1,3)$ are $2+1, \overline{2}+1$. Finally, the overpartitions in the subclass enumerated by $\overline{p}_{ou}(4)$ where the number of parts is even are $3+\overline{1},\overline{3}+\overline{1}$  whereas the set of overpartitions in the subclass enumerated by $\overline{p}_{ou}(4)$ where the number of parts is odd is empty. This gives 
$$\overline{\mathrm{spt}}1_o'(5)+\overline{\mathrm{spt}}1_o'(3)=-2=-\overline{p}_{ou}'(4).$$
\\
\section{Concluding remarks}\label{sec:secconclusions}
It appears that the degrees of the polynomials in the numerators and denominators of the rational function $\overline{P}_k(q)$ in reduced form have the same degree for a given $k$, say, $\operatorname{deg} \overline{P}_k$ and $\operatorname{deg} \overline{P}_{2k}=2k+ \operatorname{deg} \overline{P}_{2k-1}$. It would be of interest to study this further and obtain results in this spirit for the other such rational functions in Theorems \ref{thm:sptko} and \ref{thm:sptko'} as well. It would also be interesting to find combinatorial proofs of our corollaries in Section \sect{secappl}.
\\
\section{Acknowledgments} We are thankful to George Andrews for his encouragement and pointing us to his recent paper \cite{An-Ba25}. AM is partially supported by the Simons Foundation Grant TSM-00002309.

\end{document}